\documentclass[letterpaper, 10 pt, conference]{ieeeconf}  

\IEEEoverridecommandlockouts                              

\usepackage{graphics} 
\usepackage{epsfig} 
\usepackage{mathptmx} 
\usepackage{times} 
\usepackage{amsmath} 
\usepackage{amssymb}  
\usepackage{epstopdf}
\usepackage{multirow}
\usepackage{color} 

\newcommand{\notiff}{%
  \mathrel{{\ooalign{\hidewidth$\not\phantom{"}$\hidewidth\cr$\iff$}}}}

\title{\LARGE \bf
An analytical safe approximation to joint chance-constrained programming with additive Gaussian noises
}

\author{Nan Li, Ilya Kolmanovsky, and Anouck Girard
\thanks{The authors are with the Department of Aerospace Engineering,
        University of Michigan, 1320 Beal Avenue, 48109-2140 Ann Arbor, MI, USA
        {\tt\small \{nanli,ilya,anouck\}@umich.edu}}%
\thanks{This research has been supported by the National Science Foundation award CNS 1544844.}%
}
\let\labelindent\relax
\usepackage{enumitem}

\usepackage{cite}
\usepackage{graphicx}
\usepackage{amsmath}

\usepackage{color} 
\usepackage{subfig}

\begin{document}

\maketitle
\thispagestyle{empty}
\pagestyle{empty}

\begin{abstract}
We propose a safe approximation to joint chance-constrained programming where the constraint functions are additively dependent on a normally-distributed random vector. The approximation is analytical, meaning that it requires neither numerical integrations nor sampling-based probability approximations. Under mild assumptions, the approximation is a standard nonlinear program. We compare this new safe approximation to another analytical safe approximation for joint chance-constrained programming based on Boole's inequality through two examples representing the constrained control of linear Gaussian-Markov models. It is shown that our proposed safe approximation has a lower degree of conservatism compared to the one based on Boole's inequality.
\end{abstract}

\section{Introduction}

Let us consider an optimization problem of the form,
\begin{subequations}\label{equ:P01}
\begin{align}
\min_{x \in X} & \quad J(x), \\
\text{subject to} & \quad \mathbb{P}\big( F(x,\xi) \le 0 \big) \ge \beta,
\end{align}
\end{subequations}
where $x \in X \subset \mathbb{R}^{n_x}$ is the vector of optimization variables, $J: X \to \mathbb{R}$ is the cost function, $\xi$ is a random vector with probability distribution $\mathcal{P}$ supported on $\Xi \subset \mathbb{R}^{n_\xi}$, $F: X \times \Xi \to \mathbb{R}^{n_m}$ defines a set of constraints, $\mathbb{P}(A)$ denotes the probability of an event $A$, and $\beta \in [0,1]$ defines a required confidence level of constraint satisfaction. Problems in this form, introduced in \cite{charnes1958cost,miller1965chance,prekopa1970probabilistic,birge2011introduction}, are typically called {\it chance-constrained programming} problems. Chance-constrained programming has a wide range of applications, e.g., in finance \cite{dentcheva2004dual,pagnoncelli2009computational}, operational science \cite{charnes1958cost,talluri2006vendor}, and control \cite{calafiore2006scenario,mesbah2016stochastic}.

In general, problem \eqref{equ:P01} with $n_m \ge 2$, called the {\it joint chance-constrained programming} (JCCP) problem, is difficult to solve. The major difficulty lies in that evaluations of $\mathbb{P}\big( F(x,\xi) \le 0 \big)$ involve numerical integrations of multivariate distributions, which are, in general, computationally intractable.

Tractable approaches to treat JCCP problems can be classified into two groups \cite{nemirovski2012safe}: sampling-based approximations\footnote{Also in the names of scenario-based, simulation-based, Monte Carlo-based approximations.} \cite{calafiore2005uncertain,erdougan2006ambiguous,luedtke2008sample,blackmore2010probabilistic,hong2011sequential} and analytical safe approximations \cite{ben2000robust,rockafellar2000optimization,nemirovski2006convex,blackmore2009convex,grosso2014chance,paulson2017stochastic}.
The former approaches approximate a probability using a finite number of samples drawn from the distribution $\mathcal{P}$. They can be applied for an arbitrary probability distribution $\mathcal{P}$ and constraint function $F$, but have the following drawbacks: 1) They can at most provide a probabilistic guarantee of chance-constraint satisfaction as the approximation itself is random; and, 2) the required number of samples for the same level of probabilistic guarantee blows up as $\beta$ approaches one \cite{calafiore2006scenario,luedtke2008sample}. The latter approaches, analytical safe approximations, typically exploit various probability inequalities \cite{lin2011probability} to derive a deterministically constrained problem whose feasible set is contained in the feasible set of the JCCP problem \eqref{equ:P01} so that the optimal feasible solution to the new problem is a suboptimal feasible solution to \eqref{equ:P01}. Most popular choices include the exploitations of Boole's inequality and Chebyshev-Markov type inequalities \cite{nemirovski2006convex,blackmore2009convex,grosso2014chance,paulson2017stochastic}. Since inequalities are used, conservatism is introduced and the degree of conservatism substantially reflects the quality of an analytical safe approximation.

One of the most extensively investigated cases of JCCP problems is with $\mathcal{P} = \mathcal{N}$, i.e., the random vector $\xi$ follows a normal distribution, and with $F$ being additively dependent on $\xi$. Such a JCCP formulation has a broad range of applications, for instance, in the constrained control of linear Gaussian-Markov models \cite{li2000robust,li2002probabilistically,ono2008iterative,blackmore2009convex}. In this paper, we also focus on JCCP problems for such a case. We propose a new analytical safe approximation to such JCCP problems, which can be efficiently solved using standard nonlinear programming solvers without involving any numerical integrations or sampling-based probability approximations.

The notations used in this paper are standard. In particular, for a vector $\mu \in \mathbb{R}^{n}$, $\mu_j$ represents its $j$th element; for a matrix $M \in \mathbb{R}^{m \times n}$, $M_i \in \mathbb{R}^{1 \times n}$ represents its $i$th row and $M_{ij}$ represents the entry in its $i$th row, $j$th column. We use $I_n$ to represent the $n \times n$ identity matrix. Also, we use {\it iff} to represent ``if and only if'' and {\it i.i.d.} stands for ``independent and identically distributed.''

\section{Preliminaries}

In this section, we list all of the preliminary lemmas that are used to prove the main result of this paper. Although some of the lemmas listed here may be well-known, we include them for the sake of completeness and point the reader to references for their proofs.

{\bf Lemma 1 (Spectral Theorem):} Let $\Sigma \in \mathbb{R}^{n \times n}$. Then, $\Sigma^{\top} = \Sigma$ (symmetric) iff there exists $\theta \in \mathbb{R}^{n \times n}$ such that $\theta^{\top} \theta = I_n$ (orthogonal) and $\theta^{\top} \Sigma \theta = \text{diag}\big(\lambda_1, \cdots, \lambda_n \big)$ (diagonalizing $\Sigma$), where $\lambda_i \in \mathbb{R}$ for all $i = 1,\cdots,n$.

{\bf Proof:} See Theorem~7.13 of \cite{axler1997linear}. $\blacksquare$
\\

{\bf Lemma 2:} Let $\phi \sim \mathcal{N}(\mu,\Sigma)$ and set $\psi = \theta^\top \phi$, where the orthogonal matrix $\theta$ is such that $\theta^{\top} \Sigma \theta = \text{diag}\big(\lambda_1, \cdots, \lambda_n \big)$. Then, (i) $\psi \sim \mathcal{N}(\theta^{\top} \mu, \theta^{\top} \Sigma \theta)$, and (ii) the components of $\psi$ are independent.

{\bf Proof:} See Theorems~7.1 and 8.1 of \cite{gut2009probability}. $\blacksquare$
\\

{\bf Lemma 3:} Let $\big(\Omega,\mathcal{F},\mathbb{P}\big)$ be a probability space. Let $\alpha_{ij}: \Omega \to \mathbb{R}$, $i = 1,\cdots, n_1$, $j = 1,\cdots,n_2$, be random variables, and $\gamma_i \in \mathbb{R}$, $i = 1,\cdots,n_1$, be constants. Then,
\begin{equation}
\mathbb{P} \Big( \bigcap_{i = 1}^{n_1} \big(\sum_{j = 1}^{n_2} \alpha_{ij} \le \gamma_i \big)\Big) \ge \beta
\end{equation}
if there exist constants $\eta_{ij} \in \mathbb{R}$, $i = 1,\cdots, n_1$, $j = 1,\cdots,n_2$, such that
\begin{align}
& \mathbb{P}\Big( \bigcap_{i = 1}^{n_1} \big(\bigcap_{j = 1}^{n_2} (\alpha_{ij} \le \eta_{ij}) \big) \Big) \ge \beta, \nonumber \\
& \sum_{j = 1}^{n_2} \eta_{ij} \le \gamma_i, \quad i = 1,\cdots, n_1.
\end{align}

{\bf Proof:} For any $\omega \in \Omega$ such that $\alpha_{ij}(\omega) \le \eta_{ij}$ for all $i = 1,\cdots, n_1$ and $j = 1,\cdots,n_2$, it holds that $\sum_{j = 1}^{n_2} \alpha_{ij}(\omega) \le \sum_{j = 1}^{n_2} \eta_{ij} \le \gamma_i$ for all $i = 1,\cdots, n_1$. Thus, $\bigcap_{i = 1}^{n_1} \big(\bigcap_{j = 1}^{n_2} (\alpha_{ij} \le \eta_{ij}) \big) \subset \bigcap_{i = 1}^{n_1} \big(\sum_{j = 1}^{n_2} \alpha_{ij} \le \gamma_i \big)$. Therefore, $\mathbb{P} \big( \bigcap_{i = 1}^{n_1} (\sum_{j = 1}^{n_2} \alpha_{ij} \le \gamma_i) \big) \ge \mathbb{P}\big(\bigcap_{i = 1}^{n_1} (\bigcap_{j = 1}^{n_2} (\alpha_{ij} \le \eta_{ij})) \big) \ge \beta$. $\blacksquare$
\\

\section{Main result}

We consider the following JCCP problem,
\begin{subequations}\label{equ:P1}
\begin{align}
\min_x & \quad J(x), \label{equ:P11} \\
\text{subject to} & \quad \mathbb{P}\big(M \phi(x) \le m \big) \ge \beta, \label{equ:P12}
\end{align}
\end{subequations}
where $x \in \mathbb{R}^{n_x}$ is the vector of optimization variables, $J:\mathbb{R}^{n_x} \to \mathbb{R}$ is a continuously differentiable function of $x$, $\phi(x)$ is a random vector taking values in $\mathbb{R}^{n_{\phi}}$, whose distribution depends on $x$, the pair $(M,m)$, $M \in \mathbb{R}^{n_m \times n_{\phi}}$ and $m \in \mathbb{R}^{n_m}$, defines the constraint set, and $\beta \in [0,1)$ represents the required confidence level of constraint satisfaction. In particular, $\phi(x) \sim \mathcal{N}\big(\mu(x),\Sigma\big)$, i.e., $\phi(x)$ is assumed to be distributed based on a multivariate normal distribution with mean $\mu(x) \in \mathbb{R}^{n_{\phi}}$ (as a continuously differentiable function of $x$) and covariance $\Sigma \in \mathbb{R}^{n_{\phi} \times n_{\phi}}$ ($\Sigma^\top = \Sigma \succeq 0$ and independent of $x$).

Note that \eqref{equ:P12} is equivalent to $\mathbb{P}\big(M \mu(x) + M(\phi(x) - \mu(x)) \le m \big) \ge \beta$, where $\phi(x) - \mu(x) \sim \mathcal{N}\big(0,\Sigma\big)$, i.e., a zero-mean additive Gaussian noise. Note also that, without loss of generality, $n_m \ge n_\phi$, since otherwise we can redefine $\phi(x) \leftarrow M \phi(x)$ and $M \leftarrow I_{n_m}$ so that $n_m = n_\phi$. We consider the form \eqref{equ:P12} because the number of constraints, $n_m$, can be much larger than the dimension of $\phi(x)$, $n_\phi$, in many problems, and the analytical safe approximation introduced in what follows involves a set of slack variables, whose number depends only on $n_\phi$.
\newpage

{\bf Theorem 1:} Any feasible solution to the following deterministically constrained problem,
\begin{subequations}\label{equ:P2}
\begin{align}
& \min_{x,\, \{\beta_j^1,\, \beta_j^2\}_{j=1}^{n_{\phi}}} \quad J(x), \label{equ:P21} \\
& \text{subject to} \nonumber \\
& \prod_{j = 1}^{n_{\phi}} \big(\beta_{j}^1 + \beta_{j}^2 -1\big) \ge \beta, \label{equ:P22} \\
& \beta_{j}^1 + \beta_{j}^2 \ge 1, \quad j = 1, \cdots, n_{\phi}, \label{equ:P23} \\[2pt]
& 0 \le \beta_{j}^{\sigma} \le 1, \quad j = 1, \cdots, n_{\phi}, \quad \sigma = 1,2, \label{equ:P24} \\
& \sum_{j = 1}^{n_{\phi}} \Big(\sqrt{2\lambda_j}\, \big|\overline{M}_{ij}\big|\, \text{erf}^{-1} \big(2 \beta_{j}^{\sigma_{ij}} - 1\big)\Big) + \overline{M}_{i}\, \overline{\mu}(x) \nonumber \\[-1pt]
& \quad \le m_i, \quad i = 1, \cdots, n_m, \label{equ:P25}
\end{align}
\end{subequations}
is a feasible solution to the JCCP problem \eqref{equ:P1}, where $\big\{\beta_j^1, \beta_j^2\big\}_{j=1}^{n_{\phi}}$ are slack variables; $\overline{M} = M \theta$, $\overline{\mu}(x) = \theta^\top \mu(x)$, and $\lambda_j = \big(\theta^\top \Sigma \theta\big)_{jj}$, in which $\theta \in \mathbb{R}^{n_{\phi} \times n_{\phi}}$ is such that $\theta^{\top} \theta = I_{n_\phi}$ and $\theta^{\top} \Sigma \theta = \text{diag}\big(\lambda_1, \cdots, \lambda_{n_\phi} \big)$; and $\sigma_{ij} = 1$ if $\overline{M}_{ij} \ge 0$ and $\sigma_{ij} = 2$ if $\overline{M}_{ij} < 0$.
\\

{\bf Proof:} Set $\psi(x) = \theta^\top \phi(x)$. By Lemma~2(i), $\psi(x) \sim \mathcal{N}\big(\overline{\mu}(x),\overline{\Sigma} \big)$, where
\begin{equation}\label{equ:P30}
\overline{\mu}(x) = \theta^\top \mu(x), \quad \overline{\Sigma} = \theta^\top \Sigma \theta = \text{diag}\big(\lambda_1, \cdots, \lambda_{n_{\phi}}\big).
\end{equation}
By Lemma~2(ii), the components of $\psi(x)$, denoted by $\big\{\psi_1(x), \cdots, \psi_{n_\phi}(x)\big\}$, are independent.

Using $\psi(x)$, the chance constraint \eqref{equ:P12} can be written as
\begin{equation}\label{equ:P31}
\mathbb{P}\big(\overline{M} \psi(x) \le m \big) = \mathbb{P}\Big(\bigcap_{i = 1}^{n_m} \big(\sum_{j = 1}^{n_{\phi}} \overline{M}_{ij} \psi_j(x) \le m_i\big) \Big) \ge \beta,
\end{equation}
where $\overline{M} = M \theta$.

By Lemma~3, if there exists a matrix $Z \in \mathbb{R}^{n_m \times n_{\phi}}$ such that
\begin{align}
& \mathbb{P}\Big(\bigcap_{i = 1}^{n_m} \big(\bigcap_{j = 1}^{n_{\phi}} (\overline{M}_{ij} \psi_j(x) \le z_{ij}) \big) \Big) = \nonumber \\
& \mathbb{P}\Big(\bigcap_{j = 1}^{n_{\phi}} \big(\bigcap_{i = 1}^{n_m} (\overline{M}_{ij} \psi_j(x) \le z_{ij}) \big) \Big) \ge \beta, \label{equ:P32} \\
& \sum_{j = 1}^{n_{\phi}} z_{ij} \le m_i, \quad i = 1,\cdots, n_m, \label{equ:P33}
\end{align}
then \eqref{equ:P31}, and hence \eqref{equ:P12}, are satisfied.

Because $\big\{\psi_1(x), \cdots, \psi_{n_\phi}(x)\big\}$ are independent, the events $\big\{ \bigcap_{i = 1}^{n_m} (\overline{M}_{ij} \psi_j(x) \le z_{ij}) \big\}_{j = 1}^{n_\phi}$ are independent. Thus, \eqref{equ:P32} can be written as
\begin{equation}\label{equ:P34}
\prod_{j = 1}^{n_{\phi}} \mathbb{P}\Big(\bigcap_{i = 1}^{n_m} \big(\overline{M}_{ij} \psi_j(x) \le z_{ij}\big) \Big) \ge \beta,
\end{equation}
which holds iff there exists a set of probability values $\big\{\beta_1, \cdots, \beta_{n_\phi}\big\}$ such that
\begin{align}
\mathbb{P}\Big(\bigcap_{i = 1}^{n_m} \big(\overline{M}_{ij} \psi_j(x) \le z_{ij}\big) \Big) & \ge \beta_j, \quad j = 1,\cdots, n_{\phi}, \label{equ:P35}  \\
\prod_{j = 1}^{n_{\phi}} \beta_j & \ge \beta. \label{equ:P360}
\end{align}

For each $j$, we categorize $\big\{\overline{M}_{ij}\big\}_{i = 1}^{n_m}$ into three groups:
\begin{equation}
i \in \begin{cases} I_j^1 & \text{if } \overline{M}_{ij} > 0, \\
I_j^2 & \text{if } \overline{M}_{ij} < 0, \\
I_j^3 & \text{if } \overline{M}_{ij} = 0. \end{cases}
\end{equation}

Then, \eqref{equ:P35} can be written as\footnote{Here we assume $\lambda_j > 0$ to simplify the exposition. Most generally, $\lambda_j \ge 0$. The case $\lambda_j = 0$ can be considered separately, which is straightforward. It will become clear that the derivations from \eqref{equ:P41} on, and hence the final expression \eqref{equ:P460}, hold for $\lambda_j \ge 0$.}
\begin{align}
& \mathbb{P}\Big( \max_{i \in I_j^2}\frac{z_{ij}}{\overline{M}_{ij}} \le \psi_j(x) \le \min_{i \in I_j^1}\frac{z_{ij}}{\overline{M}_{ij}} \Big) = \nonumber \\
& \mathbb{P}\Big( \max_{i \in I_j^2}\frac{z_{ij} - \overline{M}_{ij} \overline{\mu}_j(x) }{\overline{M}_{ij} \sqrt{\lambda_j}}  \le \frac{\psi_j(x) - \overline{\mu}_j(x)}{\sqrt{\lambda_j}} \nonumber \\
&\quad \le \min_{i \in I_j^1}\frac{z_{ij} - \overline{M}_{ij} \overline{\mu}_j(x)}{\overline{M}_{ij} \sqrt{\lambda_j}} \Big) \ge \beta_j, \label{equ:P36} \\[3pt]
&\quad\quad z_{ij} \ge 0, \quad i \in I_j^3,
\end{align}
where $\frac{\psi_j(x) - \overline{\mu}_j(x)}{\sqrt{\lambda_j}} \sim \mathcal{N}(0,1)$.

Then, \eqref{equ:P36} can be expressed using the cumulative distribution function of the standard normal distribution, $F(\zeta) = \mathbb{P}\big(z \le \zeta\big)$, $z \sim \mathcal{N}(0,1)$, as
\begin{align}
& F \Big(\min_{i \in I_j^1}\frac{z_{ij} - \overline{M}_{ij} \overline{\mu}_j(x)}{\overline{M}_{ij} \sqrt{\lambda_j}} \Big) - F \Big( \max_{i \in I_j^2}\frac{z_{ij} - \overline{M}_{ij} \overline{\mu}_j(x) }{\overline{M}_{ij} \sqrt{\lambda_j}}\Big) \nonumber \\
&= F \Big(\min_{i \in I_j^1}\frac{z_{ij} - \overline{M}_{ij} \overline{\mu}_j(x)}{\overline{M}_{ij} \sqrt{\lambda_j}} \Big) + \nonumber \\
&\quad\, F \Big( \min_{i \in I_j^2} -\frac{z_{ij} - \overline{M}_{ij} \overline{\mu}_j(x) }{\overline{M}_{ij} \sqrt{\lambda_j}}\Big) - 1 \ge \beta_j, \label{equ:P37}
\end{align}
where we have used the property $F(\zeta) = 1 - F(-\zeta)$.

The constraint \eqref{equ:P37} holds iff there exist probability values $\beta_{j}^1$ and $\beta_{j}^2$ such that
\begin{align}
& F \Big(\min_{i \in I_j^1}\frac{z_{ij} - \overline{M}_{ij} \overline{\mu}_j(x)}{\overline{M}_{ij} \sqrt{\lambda_j}} \Big) \ge \beta_{j}^1, \label{equ:P38} \\
& F \Big( \min_{i \in I_j^2} -\frac{z_{ij} - \overline{M}_{ij} \overline{\mu}_j(x) }{\overline{M}_{ij} \sqrt{\lambda_j}}\Big) \ge \beta_{j}^2, \label{equ:P39} \\[2pt]
&\quad \beta_{j}^1 + \beta_{j}^2 -1 \ge \beta_j.
\end{align}

Using the inverse error function $\text{erf}^{-1}(\cdot)$, \eqref{equ:P38} is almost surely equivalent to
\begin{align}
& \min_{i \in I_j^1}\frac{z_{ij} - \overline{M}_{ij} \overline{\mu}_j(x)}{\overline{M}_{ij} \sqrt{\lambda_j}} \nonumber \\
& \ge F^{-1}\big(\beta_{j}^1\big) = \sqrt{2}\, \text{erf}^{-1} \big(2 \beta_{j}^1 - 1\big),
\end{align}
which is equivalent to the set of constraints
\begin{align}
& \frac{z_{ij} - \overline{M}_{ij} \overline{\mu}_j(x)}{\overline{M}_{ij} \sqrt{\lambda_j}} \ge \sqrt{2}\, \text{erf}^{-1} \big(2 \beta_{j}^1 - 1\big), \label{equ:P40} \\
& z_{ij} \ge \sqrt{2\lambda_j}\,\overline{M}_{ij}\, \text{erf}^{-1} \big(2 \beta_{j}^1 - 1\big) + \overline{M}_{ij} \overline{\mu}_j(x), \label{equ:P41}
\end{align}
for all $i \in I_j^1$, where in restating \eqref{equ:P40} as \eqref{equ:P41} we have used the fact that $\overline{M}_{ij} > 0$ for all $i \in I_j^1$.

Similarly, \eqref{equ:P39} is equivalent to the set of constraints
\begin{equation}
z_{ij} \ge -\sqrt{2\lambda_j}\,\overline{M}_{ij}\, \text{erf}^{-1} \big(2 \beta_{j}^2 - 1\big) + \overline{M}_{ij} \overline{\mu}_j(x), \label{equ:P42}
\end{equation}
for all $i \in I_j^2$. Note that $\overline{M}_{ij} < 0$ for all $i \in I_j^2$.

Combining the cases of $i \in I_j^1$, $i \in I_j^2$, and $i \in I_j^3$, we obtain
\begin{equation}\label{equ:P43}
z_{ij} \ge \sqrt{2\lambda_j}\, \big|\overline{M}_{ij}\big|\, \text{erf}^{-1} \big(2 \beta_{j}^{\sigma_{ij}} - 1\big) + \overline{M}_{ij} \overline{\mu}_j(x),
\end{equation}
for all $i = 1, \cdots, n_m$, where $\sigma_{ij} = 1$ if $\overline{M}_{ij} \ge 0$ and $\sigma_{ij} = 2$ if $\overline{M}_{ij} < 0$.

Based on \eqref{equ:P30} to \eqref{equ:P43}, we have shown that the joint chance constraint \eqref{equ:P12} is satisfied if there exist a matrix $Z \in \mathbb{R}^{n_m \times n_{\phi}}$ and a set of probability values $\big\{\beta_j, \beta_j^1, \beta_j^2\big\}_{j=1}^{n_{\phi}}$ such that
\begin{align}
& \prod_{j = 1}^{n_{\phi}} \beta_j \ge \beta, \label{equ:P441} \\
& \beta_{j}^1 + \beta_{j}^2 -1 \ge \beta_j, \quad j = 1, \cdots, n_{\phi}, \label{equ:P442} \\
& \sum_{j = 1}^{n_{\phi}} z_{ij} \le m_i, \quad i = 1, \cdots, n_m, \label{equ:P443} \\
& z_{ij} \ge \sqrt{2\lambda_j}\, \big|\overline{M}_{ij}\big|\, \text{erf}^{-1} \big(2 \beta_{j}^{\sigma_{ij}} - 1\big) + \overline{M}_{ij} \overline{\mu}_j(x), \nonumber \\[2pt]
&\quad\quad i = 1, \cdots, n_m, \quad j = 1, \cdots, n_{\phi}, \label{equ:P444}
\end{align}
where $\sigma_{ij} = 1$ if $\overline{M}_{ij} \ge 0$ and $\sigma_{ij} = 2$ if $\overline{M}_{ij} < 0$.

Furthermore, the existence of $\big\{\beta_j, \beta_j^1, \beta_j^2\big\}_{j=1}^{n_{\phi}}$ satisfying \eqref{equ:P441} and \eqref{equ:P442} is equivalent to the existence of $\big\{\beta_j^1, \beta_j^2\big\}_{j=1}^{n_{\phi}}$ satisfying
\begin{align}
& \prod_{j = 1}^{n_{\phi}} \big(\beta_{j}^1 + \beta_{j}^2 -1\big) \ge \beta, \label{equ:P450} \\
& \beta_{j}^1 + \beta_{j}^2 -1 \ge 0, \quad j = 1, \cdots, n_{\phi}, \label{equ:P45}
\end{align}
so that the variables $\big\{\beta_1, \cdots, \beta_{n_\phi}\big\}$ are dropped. In restating \eqref{equ:P441}, \eqref{equ:P442} as \eqref{equ:P450}, \eqref{equ:P45} we have used the fact that the variables $\big\{\beta_1, \cdots, \beta_{n_\phi}\big\}$ are probability values, thus, must take non-negative values.

Similarly, \eqref{equ:P443} and \eqref{equ:P444} are equivalent to
\begin{align}
& \sum_{j = 1}^{n_{\phi}} \Big(\sqrt{2\lambda_j}\, \big|\overline{M}_{ij}\big|\, \text{erf}^{-1} \big(2 \beta_{j}^{\sigma_{ij}} - 1\big) + \overline{M}_{ij} \overline{\mu}_j(x)\Big) \nonumber \\
& \quad\quad \le m_i, \quad i = 1, \cdots, n_m,
\end{align}
where $\sigma_{ij} = 1$ if $\overline{M}_{ij} \ge 0$ and $\sigma_{ij} = 2$ if $\overline{M}_{ij} < 0$, so that the variables $Z \in \mathbb{R}^{n_m \times n_{\phi}}$ are dropped.

To sum up, if $\big(x,\{\beta_j^1, \beta_j^2\}_{j=1}^{n_{\phi}}\big)$ satisfies the following set of deterministic constraints,
\begin{subequations}\label{equ:P460}
\begin{align}
& \prod_{j = 1}^{n_{\phi}} \big(\beta_{j}^1 + \beta_{j}^2 -1\big) \ge \beta, \\
& \beta_{j}^1 + \beta_{j}^2 \ge 1, \quad j = 1, \cdots, n_{\phi}, \\[2pt]
& 0 \le \beta_{j}^{\sigma} \le 1, \quad j = 1, \cdots, n_{\phi}, \quad \sigma = 1,2, \label{equ:P46} \\
& \sum_{j = 1}^{n_{\phi}} \Big(\sqrt{2\lambda_j}\, \big|\overline{M}_{ij}\big|\, \text{erf}^{-1} \big(2 \beta_{j}^{\sigma_{ij}} - 1\big)\Big) + \overline{M}_{i}\, \overline{\mu}(x) \nonumber \\[-1pt]
& \quad \le m_i, \quad i = 1, \cdots, n_m,
\end{align}
\end{subequations}
where $\sigma_{ij} = 1$ if $\overline{M}_{ij} \ge 0$ and $\sigma_{ij} = 2$ if $\overline{M}_{ij} < 0$, then $x$ satisfies the joint chance constraint
\begin{equation}
\mathbb{P}\big(M \phi(x) \le m \big) \ge \beta.
\end{equation}
Note that the constraints \eqref{equ:P46} come from the fact that the variables $\big\{\beta_j^1, \beta_j^2\big\}_{j=1}^{n_{\phi}}$ are probability values. This completes the proof.
$\blacksquare$
\\

Note that $\theta$, $\overline{M}$, $\big\{\lambda_j\big\}_{j=1}^{n_{\phi}}$, and the set of indices $\big\{\sigma_{ij}\big\}_{i = 1,\cdots,n_m,\, j = 1, \cdots, n_{\phi}}$ are independent of $x \in \mathbb{R}^{n_x}$ and $\big\{\beta_j^1, \beta_j^2\big\}_{j=1}^{n_{\phi}}$, and thus can be determined before solving the optimization problem \eqref{equ:P2}.

The significance of Theorem~1 is as follows. On the one hand, the deterministically constrained problem \eqref{equ:P2} is an analytical safe approximation to the original joint chance-constrained problem \eqref{equ:P1}. On the other hand, \eqref{equ:P2} is a standard nonlinear programming problem with continuously differentiable cost and constraint functions, which does not require numerical integrations or sampling-based probability approximations and can hence be solved using standard nonlinear programming solvers.

We note that the derivation of \eqref{equ:P2} relies on the diagonalization of a multivariate normal distribution in \eqref{equ:P30} and the exploitation of independent events in \eqref{equ:P32} $\implies$ \eqref{equ:P34}. An alternative approach could be considered, which is to directly diagonalize $M \phi(x)$ as $M \phi(x) = \theta' \psi'(x)$ with $\psi'(x) \sim \mathcal{N}\big(\mu'(x),\text{diag}(\lambda_1', \cdots, \lambda_{n_m}')\big)$ so that \eqref{equ:P12} can be written as $\mathbb{P}\big(\theta' \psi'(x) \le m \big) \ge \beta$. One may hope to work with $\mathbb{P}\big(\psi'(x) \le (\theta')^{-1} m \big) = \prod_{j=1}^{n_m} \mathbb{P}\big(\psi_j'(x) \le ((\theta')^{-1} m)_j \big) \ge \beta$. However, in general $\theta' \psi'(x) \le m \notiff \psi'(x) \le (\theta')^{-1} m$ even if $\theta'$ is orthogonal. A straightforward counterexample is $\begin{bmatrix} -\frac{\sqrt{3}}{2} & \frac{1}{2} \\ \frac{1}{2} & \frac{\sqrt{3}}{2} \end{bmatrix} \begin{bmatrix} 1 \\ -1 \end{bmatrix} \le \begin{bmatrix} 0 \\ 0 \end{bmatrix}$.
\\

\section{Risk allocation using Boole's inequality}

Two alternative approaches to obtain analytical safe approximations to the joint chance-constrained problem \eqref{equ:P1} that have been extensively studied in the literature are: 1) to separate the joint constraint into multiple elementary constraints and exploit Boole's inequality \cite{ono2008iterative,blackmore2009convex,paulson2017stochastic}, and 2) to ensure that the $\beta$-level confidence ellipsoid of the random vector $\phi(x)$ is contained in the constraint admissible set \cite{van2001lmi,van2006stochastic}. It is possible to show that, in general, the second method is more conservative than the first method \cite{blackmore2009convex,cinquemani2011convexity}. Thus, we choose to use the first method as a comparison benchmark. It is briefly reviewed in this section.

The constraint \eqref{equ:P12} can be written as
\begin{align}
& \mathbb{P}\Big(\bigcap_{i = 1}^{n_m} \big(M_{i}\, \phi(x) \le m_i\big) \Big) \ge \beta, \\
& \mathbb{P}\Big(\bigcup_{i = 1}^{n_m} \big(M_{i}\, \phi(x) > m_i\big) \Big) \le 1 - \beta,
\end{align}
which holds if there exists a set of probability values $\big\{\beta_1, \cdots, \beta_{n_m}\big\}$ such that
\begin{align}
\mathbb{P}\big(M_{i}\, \phi(x) \le m_i\big) &= \mathbb{P}\Big(\frac{M_{i}\, \phi(x) - M_{i}\, \mu(x)}{\sqrt{M_{i}\, \Sigma M_{i}^\top}} \le \frac{m_i - M_{i}\, \mu(x)}{\sqrt{M_{i}\, \Sigma M_{i}^\top}}\Big) \nonumber \\[-1pt]
&\ge \beta_i, \quad i = 1,\cdots,n_m, \label{equ:P51} \\[2pt]
& \sum_{i = 1}^{n_m} \big(1 - \beta_i\big) \le 1 - \beta,
\end{align}
where $\frac{M_{i}\, \phi(x) - M_{i}\, \mu(x)}{\sqrt{M_{i}\, \Sigma M_{i}^\top}} \sim \mathcal{N}(0,1)$, which is based on Boole's inequality:
\begin{align}
& \mathbb{P}\Big(\bigcup_{i = 1}^{n_m} \big(M_{i}\, \phi(x) > m_i\big) \Big) \le \sum_{i = 1}^{n_m} \Big(\mathbb{P}\big(M_{i}\, \phi(x) > m_i\big) \Big).
\end{align}

Using the inverse error function $\text{erf}^{-1}(\cdot)$, \eqref{equ:P51} is almost surely equivalent to
\begin{align}
& \frac{m_i - M_{i}\, \mu(x)}{\sqrt{M_{i}\, \Sigma M_{i}^\top}} \ge F^{-1}\big(\beta_i\big) = \sqrt{2}\, \text{erf}^{-1} \big(2 \beta_i - 1\big), \nonumber \\
& M_{i}\, \mu(x) \le m_i - \sqrt{2 M_{i}\, \Sigma M_{i}^\top}\, \text{erf}^{-1} \big(2 \beta_i - 1\big).
\end{align}

To sum up, any feasible solution to the following deterministically constrained problem,
\begin{subequations}\label{equ:P5}
\begin{align}
& \min_{x,\, \{\beta_i\}_{i=1}^{n_m}} \quad J(x), \\
& \text{subject to} \nonumber \\
& \sum_{i = 1}^{n_m} \big(1 - \beta_i\big) \le 1 - \beta, \\
& 0 \le \beta_i \le 1, \quad i = 1, \cdots, n_m, \\[2pt]
& M_{i}\, \mu(x) \le m_i - \sqrt{2 M_{i}\, \Sigma M_{i}^\top}\, \text{erf}^{-1} \big(2 \beta_i - 1\big), \nonumber \\[2pt]
&\quad\quad\quad\quad i = 1, \cdots, n_m,
\end{align}
\end{subequations}
is a feasible solution to the JCCP problem \eqref{equ:P1}, where $\{\beta_i\}_{i=1}^{n_m}$ are slack variables.

Note that, similar to \eqref{equ:P2}, \eqref{equ:P5} is also a nonlinear programming problem with continuously differentiable cost and constraint functions.

\section{Application to constrained control of linear Gaussian-Markov models}

In this section, we use examples representing the constrained control of linear Gaussian-Markov models to illustrate the effectiveness of the analytical safe approximation \eqref{equ:P2} to the joint chance-constrained programming problem \eqref{equ:P1}.

Consider a discrete-time linear Gaussian-Markov model,
\begin{subequations}\label{equ:P6}
\begin{align}
x_{t+1} &= A x_t + B_u u_t + B_w w_t, \label{equ:P61} \\[1pt]
y_t &= C x_t + D_u u_t + D_w w_t, \\[1pt]
x_0 &\sim \mathcal{N}(\overline{x}_0,\Sigma_x), \\[1pt]
w_t &\sim \mathcal{N}(0,\Sigma_w), \quad t \in \mathbb{Z}_{\ge 0}, \label{equ:P64}
\end{align}
\end{subequations}
where 1) $\overline{x}_0$, $\Sigma_x$, and $\Sigma_w$ are given, 2) $\{w_t\}_{t \in \mathbb{Z}_{\ge 0}}$ are {\it i.i.d.}, and 3) $\{w_t\}_{t \in \mathbb{Z}_{\ge 0}}$ are independent of $x_0$.

The control objective is to minimize a quadratic cost function,
\begin{align}\label{equ:P91}
& J(u_0,\cdots,u_{N-1}) = \mathbb{E} \Big[ \sum_{t = 0}^{N-1} \big(x_{t+1}^\top Q x_{t+1} + u_t^\top R u_t\big) \Big], \nonumber \\[-1pt]
& = \sum_{t = 0}^{N-1} \big(\mathbb{E}[x_{t+1}]^\top Q \mathbb{E}[x_{t+1}] + u_t^\top R u_t\big) + const.,
\end{align}
where $Q^\top = Q \succeq 0$, $R^\top = R \succ 0$, $N$ is the prediction horizon, and $const.$ represents constant terms independent of $\{u_0,\cdots,u_{N-1}\}$, subject to the joint chance constraint,
\begin{equation}\label{equ:P92}
\mathbb{P}\Big( \bigcap_{t = 1}^{N} \big(y_t \le y_{\max} \big) \Big) \ge \beta.
\end{equation}

We let $\phi = \begin{bmatrix} y_1^\top & \cdots & y_N^\top \end{bmatrix}^\top$. The covariance of $y_i$ and $y_j$, $i \le j$, is given by
\small
\begin{align}
& \Sigma(i,j) = \mathbb{E}\Big(\big(y_i - \mathbb{E}[y_i]\big)\big(y_j - \mathbb{E}[y_j]\big)^\top\Big) \\
=&\, \mathbb{E}\bigg( \begin{bmatrix} C A^{i} \mkern-8mu & C A^{i-1} B_w \mkern-8mu & \cdots \mkern-8mu & CB_w \mkern-8mu & D_w \end{bmatrix} \begin{bmatrix} x_0 - \overline{x}_0 \\ w_{0} \\ \vdots \\ w_{i-1} \\  w_i \end{bmatrix} \nonumber \\[-3pt]
&\quad \begin{bmatrix} x_0 - \overline{x}_0 \\ w_{0} \\ \vdots \\ w_{j-1} \\  w_j \end{bmatrix}^\top \begin{bmatrix} C A^{j} \mkern-8mu & C A^{j-1} B_w \mkern-8mu & \cdots \mkern-8mu & CB_w \mkern-8mu & D_w \end{bmatrix}^\top\bigg) \\[3pt]
=&\, \begin{bmatrix} C A^{i} \mkern-8mu & C A^{i-1} B_w \mkern-8mu & \cdots \mkern-8mu & CB_w \mkern-8mu & D_w \end{bmatrix} \nonumber \\
&\begin{bmatrix} \,\,\Sigma_x & & \\ & \begin{bmatrix} \Sigma_w & \\ & \ddots & \\ & & \Sigma_w \end{bmatrix} & 0\,\, \end{bmatrix} \begin{bmatrix} C A^{j} \mkern-8mu & C A^{j-1} B_w \mkern-8mu & \cdots \mkern-8mu & CB_w \mkern-8mu & D_w \end{bmatrix}^\top,
\end{align}
\normalsize
which is used to construct the covariance matrix $\Sigma$ in the JCCP formulation \eqref{equ:P1}. Once $\Sigma$ is obtained, we solve a standard eigenvalue problem to obtain the eigenvalues $\{\lambda_1, \cdots, \lambda_n\}$ and a set of corresponding orthonormal eigenvectors $\{\nu_1, \cdots, \nu_n\}$ of $\Sigma$, after which the orthogonal matrix $\theta$ is constructed by $\theta = \begin{bmatrix} \nu_1, \cdots, \nu_n \end{bmatrix}$.

After transforming the problem \eqref{equ:P91} and \eqref{equ:P92} into the form of \eqref{equ:P1}, and further into its analytical safe approximation in the form of \eqref{equ:P2}, we use the standard nonlinear programming solver Matlab $fmincon$ function with the interior-point method \cite{boyd2004convex} to solve for $\{u_0,\cdots,u_{N-1}\}$.

To evaluate the effectiveness of our analytical safe approximation \eqref{equ:P2}, once an optimal solution $\{u_0,\cdots,u_{N-1}\}$ is obtained, we apply it to the open-loop system \eqref{equ:P6} where the disturbance signals $\{w_0,\cdots,w_{N-1}\}$ are randomly created based on \eqref{equ:P64}, and repeat such a simulation for $10,000$ times.

We note that here we consider open-loop control. The use of the approach in the setting of receding-horizon optimal control, e.g., in stochastic model predictive control \cite{mesbah2016stochastic}, to achieve closed-loop operation represents a natural extension, which is left as a topic to future research.

To compare the performance of our analytical safe approximation \eqref{equ:P2} and the one  \eqref{equ:P5} based on the exploitation of Boole's inequality, we also use \eqref{equ:P5} to solve for $\{u_0,\cdots,u_{N-1}\}$ and run the same experiment.

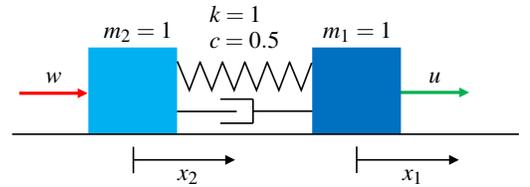
\begin{figure}[!ht]
\begin{center}
\begin{picture}(180, 66.0)
\put(0,0){\epsfig{file=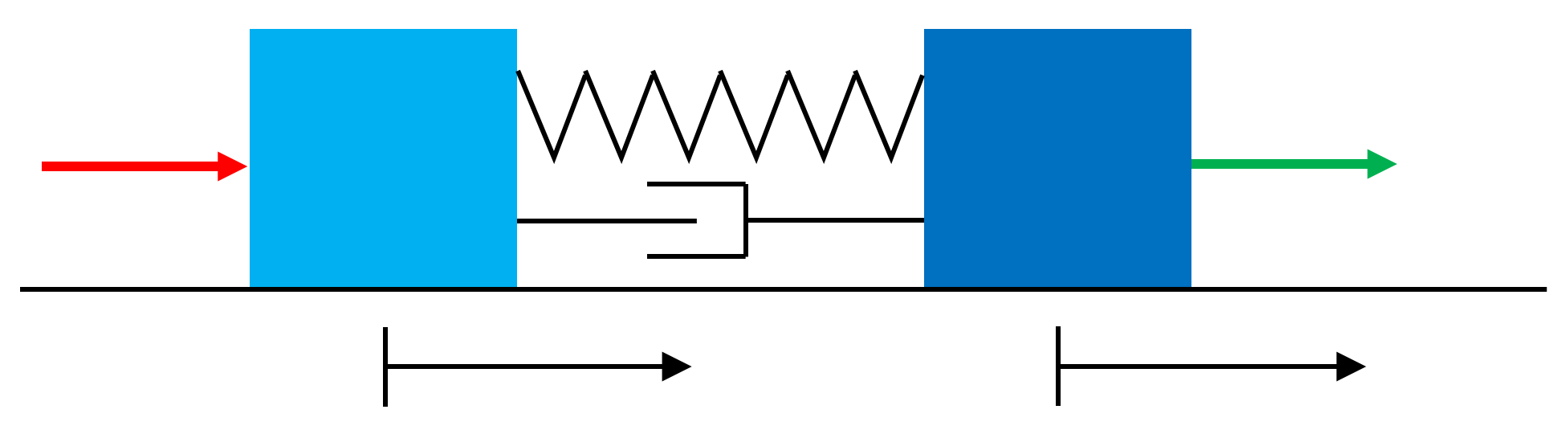, width = 0.8\linewidth}}
\small
\put(150,0){$x_1$}
\put(65,0){$x_2$}
\put(77,60){$k = 1$}
\put(77,50){$c = 0.5$}
\put(120,54){$m_1 = 1$}
\put(37,54){$m_2 = 1$}
\put(160,37){$u$}
\put(15,37){$w$}
\normalsize
\end{picture}
\end{center}
      \caption{{\small Double mass-spring-damper system.}}
      \label{fig:mass-spring-damper}
\end{figure}

Comparison results are based on the following two examples:

{\it Example 1:} We consider the following model representing the double mass-spring-damper system shown in Fig.~\ref{fig:mass-spring-damper},
\begin{align}\label{equ:EX1}
    \frac{d}{dt} \begin{bmatrix} x_1 \\ x_2 \\ \dot{x}_1 \\ \dot{x}_2 \end{bmatrix} =&\, \begin{bmatrix} 0 & 0 & 1 & 0 \\ 0 & 0 & 0 & 1 \\ -\frac{k}{m_1} & \frac{k}{m_1} & -\frac{c}{m_1} & \frac{c}{m_1} \\ \frac{k}{m_2} & -\frac{k}{m_2} & \frac{c}{m_2} & -\frac{c}{m_2}\end{bmatrix} \begin{bmatrix} x_1 \\ x_2 \\ \dot{x}_1 \\ \dot{x}_2 \end{bmatrix} \nonumber \\
    &\,\, + \begin{bmatrix} 0 & 0 \\ 0 & 0 \\ \frac{1}{m_1} & 0 \\ 0 & \frac{1}{m_2} \end{bmatrix} \begin{bmatrix} u \\ w \end{bmatrix}.
\end{align}
We discretize \eqref{equ:EX1} using the Matlab $c2d$ function with sample time of $\Delta t = 0.5$~[sec] to obtain the corresponding discrete-time model in the form of \eqref{equ:P61}. We consider
\begin{align}\label{equ:EX1_2}
& C = \begin{bmatrix} 1 & 0 & 0 & 0 \\
    0 & 1 & 0 & 0 \end{bmatrix}, \quad D_u = D_w = 0, \\
& \overline{x}_0 = \begin{bmatrix} -0.5 & -0.5 & 0 & 0 \end{bmatrix}^\top, \quad \Sigma_x = 0, \quad \Sigma_w = 10^{-4}, \nonumber
\end{align}
and
\begin{align}
& Q = \text{diag}(1000,1000,1,1), \quad\quad R = 1, \nonumber \\[0.5pt]
& N = 20, \quad\quad y_{\max} = \begin{bmatrix} 0 & 0 \end{bmatrix}^\top.
\end{align}

\begin{figure}[h!]
\begin{center}
\begin{picture}(246.0, 230.0)
\put(  0,  116){\epsfig{file=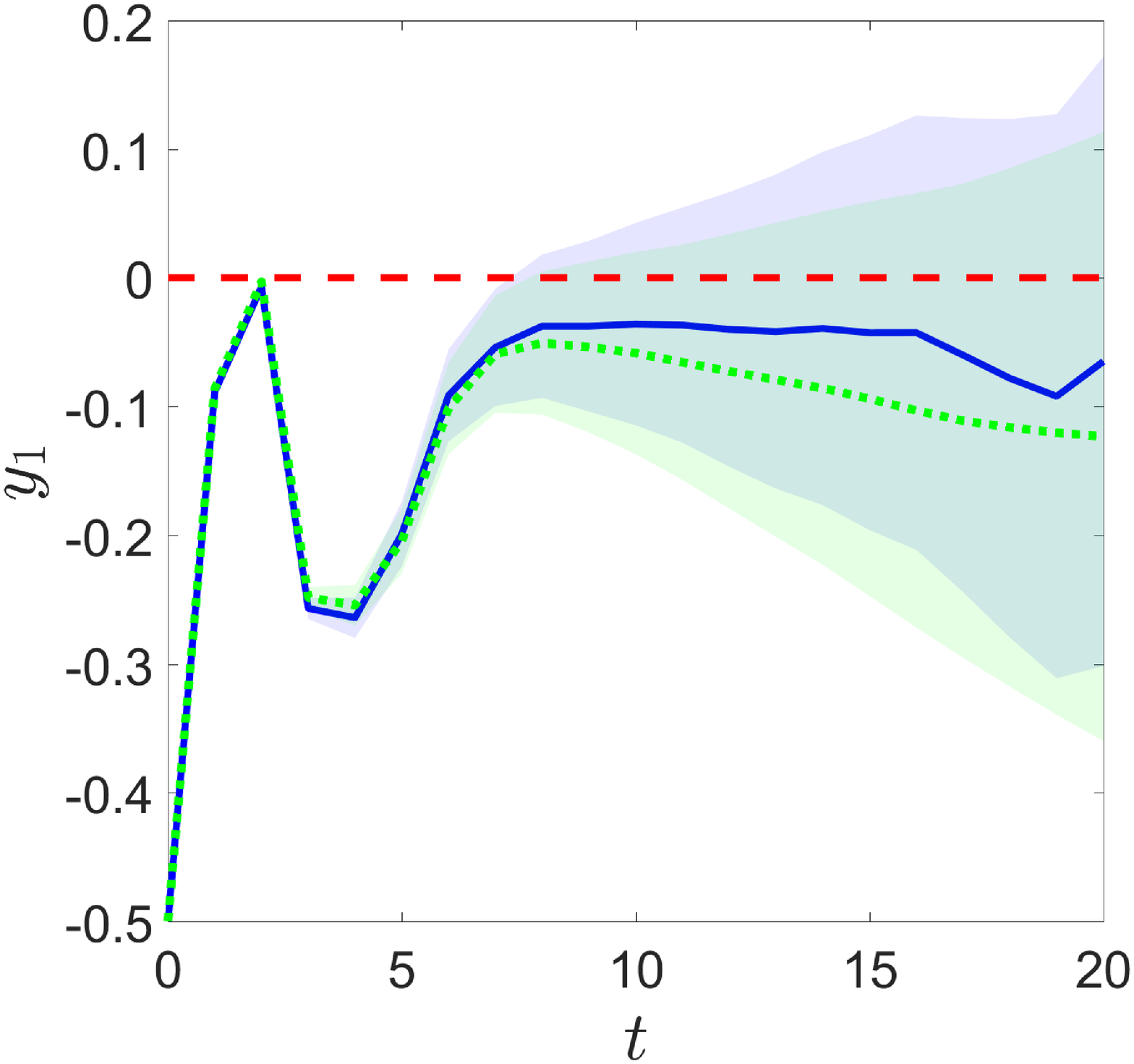,height=1.62in}}  
\put(  123,  116){\epsfig{file=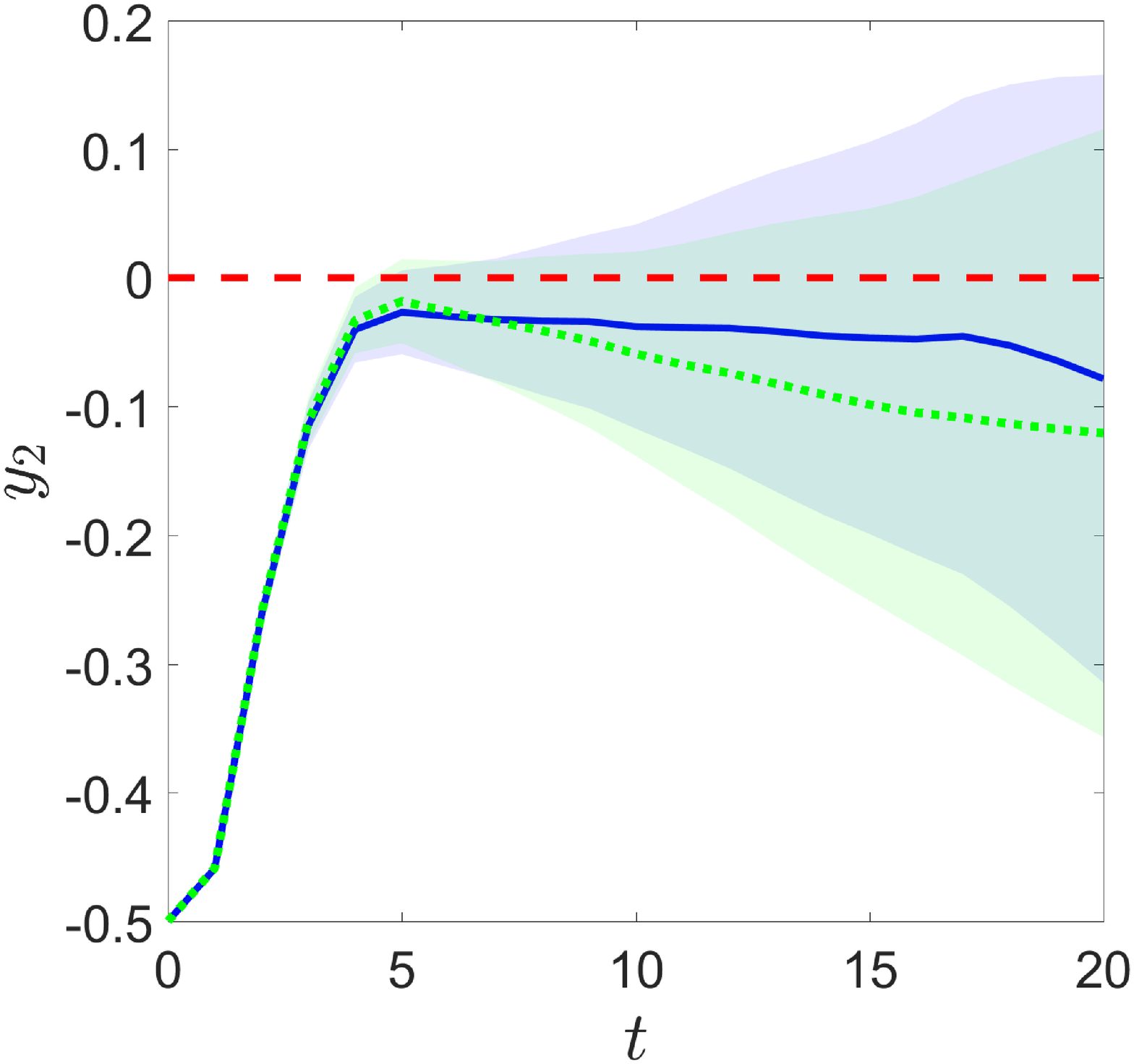,height=1.62in}}  
\put(  0,  0){\epsfig{file=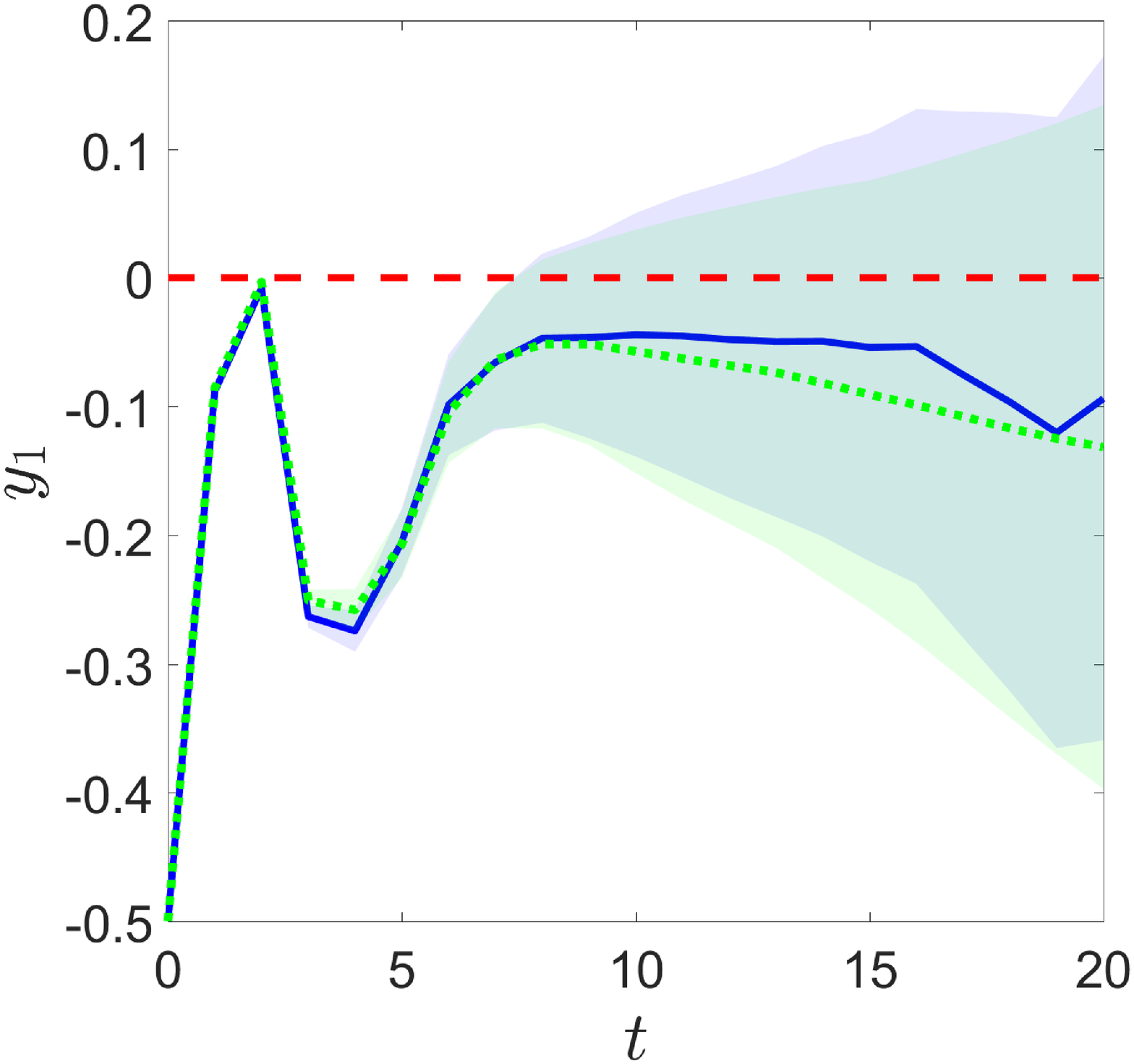,height=1.62in}}  
\put(  123,  0){\epsfig{file=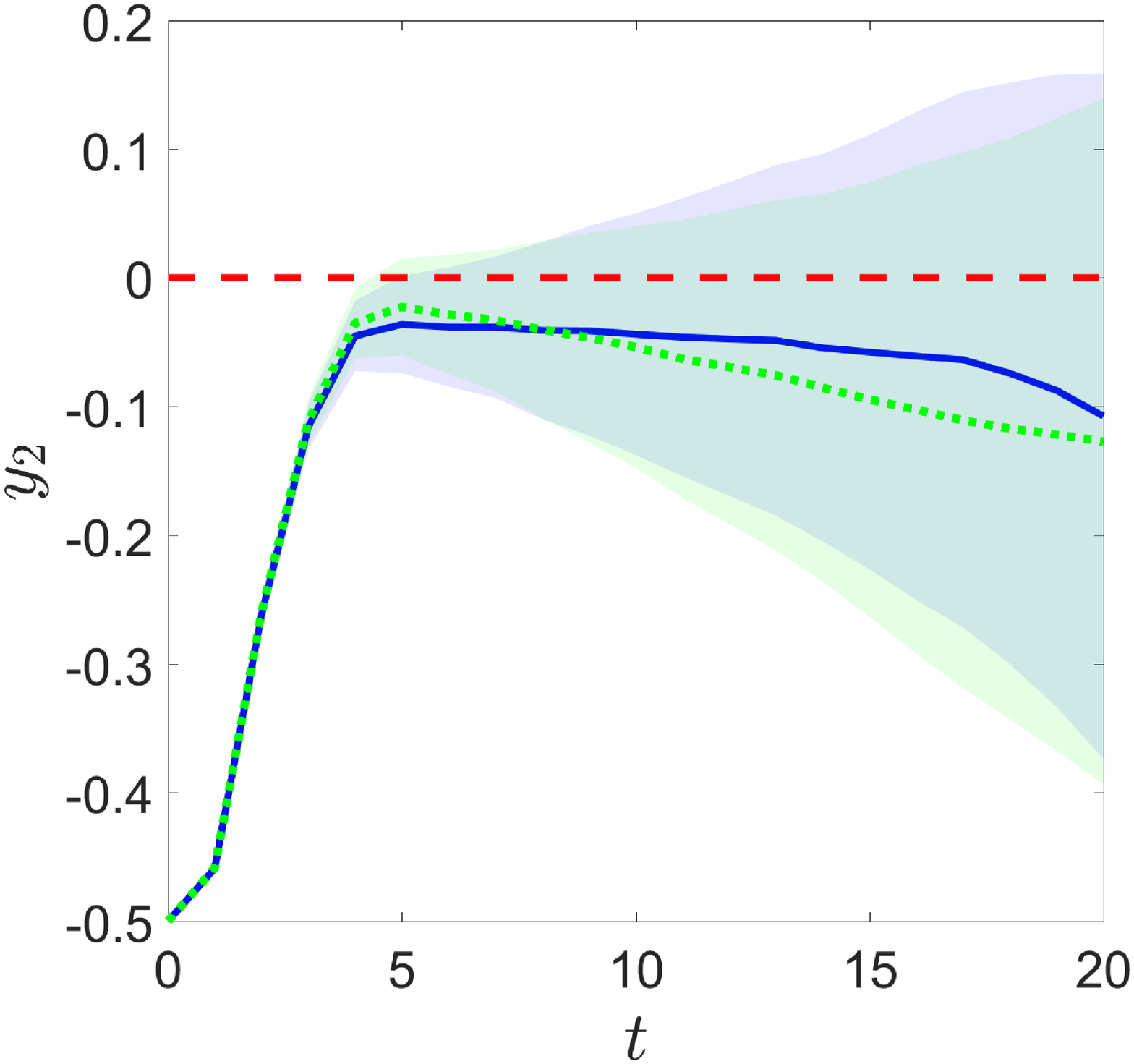,height=1.62in}}  
\small
\put(  225,  226){(a)}  
\put(  225,  110){(b)}  
\normalsize
\end{picture}
\end{center}
      \caption{{\small Simulation results of {\it Example~1} with (a) $\beta = 0.6$ and (b) $\beta = 0.8$. The blue-solid (resp. green-dotted) curves represent the expected trajectories of $y$ under the control input sequences $\{u_0,\cdots,u_{N-1}\}$ solved based on \eqref{equ:P2} (resp. \eqref{equ:P5}). The light blue (resp. green) shadowed areas represent their corresponding envelopes of $10,000$ measured trajectories. The red-dashed lines represent the constraint boundaries $y_{\max}$.}}
      \label{fig:1}
\end{figure}

We test two cases: $\beta = 0.6$ and $\beta = 0.8$. The responses of $y$ under the control input sequences $\{u_0,\cdots,u_{N-1}\}$ solved based on \eqref{equ:P2} and \eqref{equ:P5} are shown in Fig.~\ref{fig:1}. It can be observed that, for both cases, $y_1$ and $y_2$ corresponding to the solutions of \eqref{equ:P2} get closer to the constraint boundaries compared to those corresponding to the solutions of \eqref{equ:P5}.

We use two metrics to compare the relative degree of conservatism between \eqref{equ:P2} and \eqref{equ:P5}. They are the cost values $J$ and the measured rates of constraint satisfaction $\overline{\beta}$, i.e., the proportion of simulation runs where the constraint $\bigcap_{t = 1}^{N} (y_t \le y_{\max})$ is satisfied, corresponding to the solutions of \eqref{equ:P2} and \eqref{equ:P5}.
Note that since both solutions are feasible solutions to the original problem \eqref{equ:P91} and \eqref{equ:P92}, their corresponding cost values reflect their relative degree of conservatism.
The comparison results for both cases are summarized in Table~\ref{tab:1}. It can be observed that, for both cases, the solution of \eqref{equ:P2} has a lower cost value and a measured rate of constraint satisfaction closer to the required value $\beta$ compared to the solution of \eqref{equ:P5}.

\begin{table}[ht]
    \centering
    \caption{{\small Comparison results of {\it Example~1}.}}
    \begin{tabular}{|c|c|c|c|}
        \hline
        \multicolumn{1}{|c}{} & & Sol. of \eqref{equ:P2} & Sol. of \eqref{equ:P5} \\
        \hline
        \multirow{2}{*}{$\beta = 0.6$} & $J$ & 597.7 & 729.7 \\ \cline{2-4}
        & $\overline{\beta}$ & 0.7737 & 0.9577 \\
        \hline
        \multirow{2}{*}{$\beta = 0.8$} & $J$ & 695.9 & 788.4 \\ \cline{2-4}
        & $\overline{\beta}$ & 0.9107 & 0.9782 \\
        \hline
    \end{tabular}
    \label{tab:1}
\end{table}

{\it Example 2:} We consider the following model,
\begin{align} 
& A = \begin{bmatrix} 1.0000  &  0.1025  &  0.2080  &  -0.0502  & -0.0057 \\
         0  &  1.1175   &  4.1534  & -0.8000  & -0.1010 \\
         0  &  0.0955   &  1.0722  & -0.0541  & -0.0153 \\
         0   &      0     &    0   &  0.1353    &     0 \\
         0   &      0     &    0    &     0   &  0.1353 \end{bmatrix}, \nonumber
\end{align}
\begin{align}\label{equ:P7}
& B_u = B_w = \begin{bmatrix} -0.0377 &  -0.0040 \\
   -1.0042  &  -0.1131 \\
   -0.0453  &  -0.0175 \\
    0.8647   &      0 \\
         0  &   0.8647 \end{bmatrix},
\end{align}
which describes the short-period pitch attitude dynamics augmented by control actuator dynamics (elevator and flaperons) of an AFTI/F-16 aircraft at the flight condition of altitude $3000$~[feet] and Mach number $0.6$. A continuous-time model is taken from \cite{sobel1985design} and discretized with sample time of $\Delta t=0.1$~[sec] to obtain the discrete-time model \eqref{equ:P7}. Also, we consider
\begin{align}\label{equ:P8}
& C = \begin{bmatrix} -1 & 0 & 0 & 0 & 0 \\
    0 & -1 & 0 & 0 & 0 \end{bmatrix}, \quad D_u = D_w = 0, \\
& \overline{x}_0 = \begin{bmatrix} 1 & 0 & 0 & 0 & 0 \end{bmatrix}^\top, \quad \Sigma_x = 0, \quad \Sigma_w = 2.5 \times 10^{-3} I_2. \nonumber
\end{align}
and
\begin{align}
& Q = \text{diag}(1000,1,1,1,1), \quad R = \text{diag}(1,1), \nonumber \\[0.5pt]
& N = 10, \quad\quad y_{\max} = \begin{bmatrix} 0 & 1 \end{bmatrix}^\top.
\end{align}

\begin{figure}[h!]
\begin{center}
\begin{picture}(246.0, 123.0)
\put(  0,  0){\epsfig{file=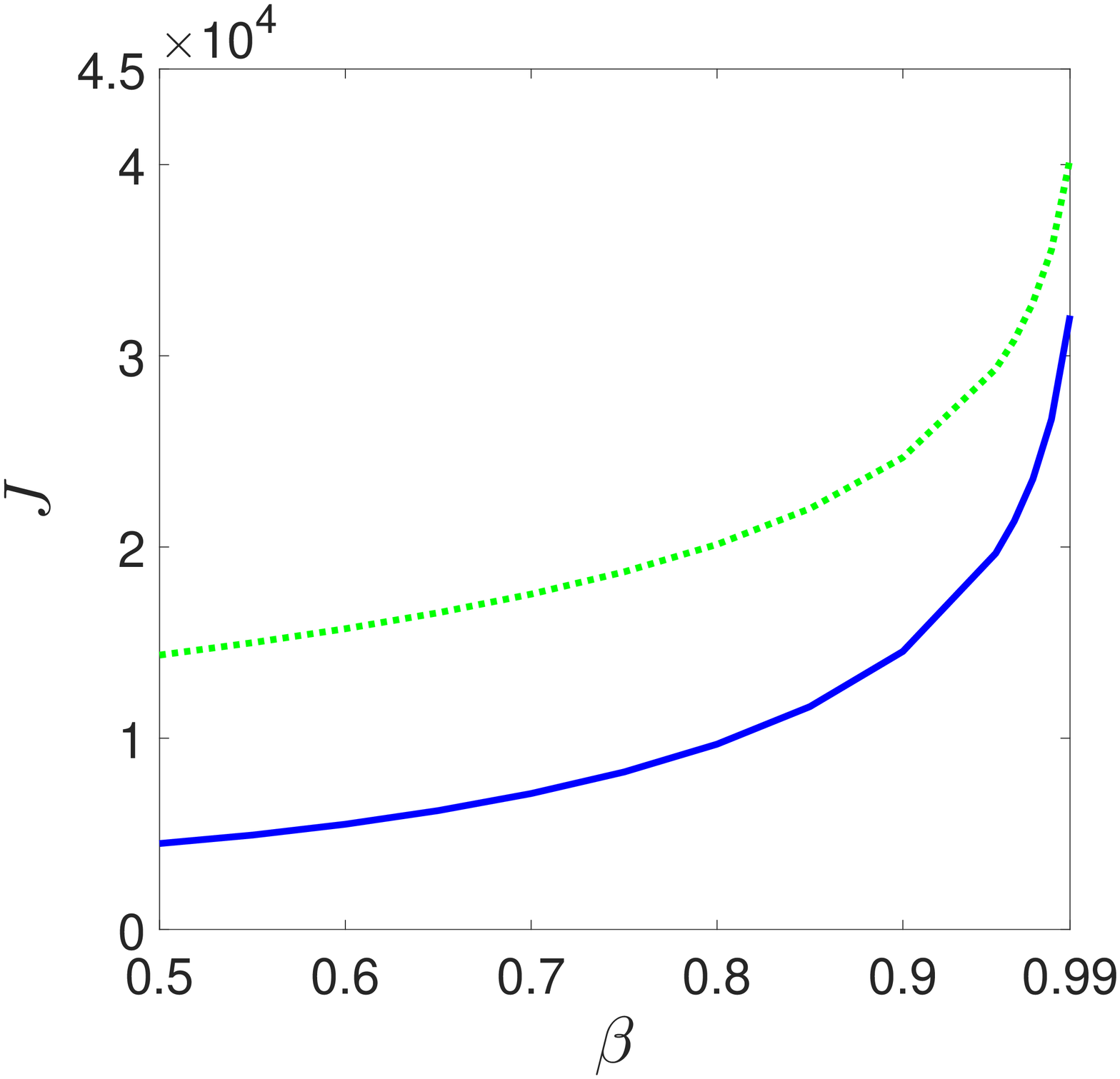,height=1.62in}}  
\put(  123,  0){\epsfig{file=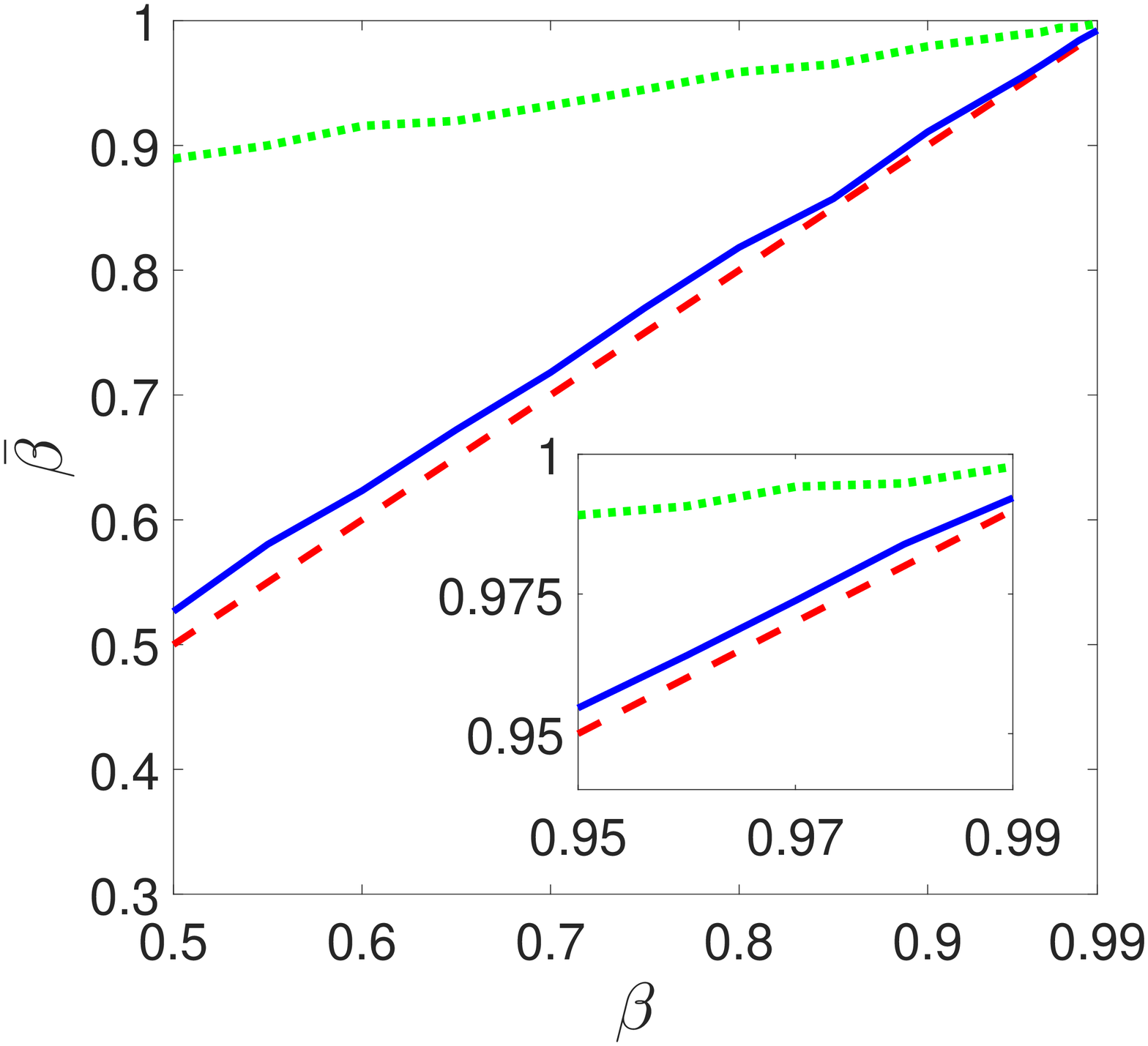,height=1.62in}}  
\small
\put(  102,  112){(a)}  
\put(  225,  112){(b)}  
\normalsize
\end{picture}
\end{center}
      \caption{{\small Comparison results of {\it Example~2}. (a) Cost values $J$ and (b) measured rates of constraint satisfaction $\overline{\beta}$ for different values of required confidence level $\beta$. The blue-solid curves correspond to the solutions of \eqref{equ:P2}, the green-dotted curves correspond to the solutions of \eqref{equ:P5}, and the red-dashed line in (b) represents the line of $\overline{\beta} = \beta$.}}
      \label{fig:2}
\end{figure}

Similar to {\it Example~1}, we use the cost values $J$ and the measured rates of constraint satisfaction $\overline{\beta}$ corresponding to the solutions of \eqref{equ:P2} and \eqref{equ:P5} for different values of required confidence level $\beta \in [0.5,0.99]$ to compare their relative degree of conservatism. The comparison results are plotted in Fig.~\ref{fig:2}.

It can be observed that the solutions of \eqref{equ:P2} have significantly lower cost values than the solutions of \eqref{equ:P5}. This can be explained by the fact that the measured rates of constraint satisfaction $\overline{\beta}$ corresponding to the solutions of \eqref{equ:P2} are much closer to the required values $\beta$ than those corresponding to the solutions of \eqref{equ:P5}. As the solutions of \eqref{equ:P2} are capable of getting closer to the constraint boundaries, they lead to lower cost values.

Therefore, the above two examples illustrate that our analytical safe approximation \eqref{equ:P2} can have a considerably lower degree of conservatism than the analytical safe approximation \eqref{equ:P5} based on the exploitation of Boole's inequality.

\section{Discussions}

There are two challenges faced by our analytical safe approximation \eqref{equ:P2} compared to the analytical safe approximation \eqref{equ:P5} based on the exploitation of Boole's inequality: 1) Under additional assumptions, including that the cost function is convex in $x$, the mean $\mu(x)$ is affine in $x$, and the required confidence level satisfies $\beta \ge 0.5$, it is possible to show that \eqref{equ:P5} is convex \cite{blackmore2009convex,paulson2017stochastic}; on the contrary, \eqref{equ:P2} is not convex in general due to the constraint \eqref{equ:P22}. And, 2) the formulation \eqref{equ:P2} involves more slack variables than \eqref{equ:P5} when $2 n_\phi > n_m$. Both 1) and 2) result in higher computational complexity of \eqref{equ:P2} than \eqref{equ:P5}. However, the computational effort in solving \eqref{equ:P2} can usually be manageable and be much less than that in solving the original JCCP problem \eqref{equ:P1} using an approach based on numerical or sampling-based integrations of multivariate distributions. For instance, the average and worst computation times for solving \eqref{equ:P2} of {\it Example 2} are $4.933$~[sec] and $5.515$~[sec] using the Matlab $fmincon$ function with the interior-point method in uncompiled code on a PC with Intel Core i7-4790 3.60 GHz processor and 16.0 GB RAM. For comparison, the average and worst computation times for solving \eqref{equ:P5} of {\it Example 2} in the same computation environment are $1.766$~[sec] and $2.666$~[sec].

\section{Conclusions}

In this paper, we proposed a new analytical safe approximation to joint chance-constrained programming problems with constraint functions additively dependent on normally-distributed random vectors. The approximation is a standard nonlinear program with continuously differentiable cost and constraint functions. Two examples representing the constrained control of linear Gaussian-Markov models were used to illustrate the effective application of our proposed analytical safe approximation and that our approximation can have a considerably lower degree of conservatism compared to a popularly used analytical safe approximation based on the exploitation of Boole's inequality.

\section*{}
\bibliographystyle{IEEEtran}
\bibliography{ref}

\end{document}